\documentclass{article}
\usepackage{amsmath}
\usepackage{amssymb}
\usepackage{amsthm}
\newtheorem{thm}{Theorem}

\begin{document}

\title{A Correction Note: Attractive Nearest Neighbor Spin Systems on the Integers}

\author{Jeffrey Lin (UCLA)}

\maketitle

\begin{abstract}
In this note, we discuss a proof in T. Liggett's work on attractive translation invariant nearest neighbor spin systems on the integers.  A correction to a wrong estimate is provided.
\end{abstract}

\section{Introduction}
For the purposes of this note, an attractive translation invariant nearest neighbor spin system is a certain kind of Feller process defined on the compact state space $X=\{0, 1\}^\mathbb{Z}$ with rates satisfying the attractiveness inequalities, and depending only on the nearest neighbors. (See \cite{book}, p.144-145 for more details.)

T. Liggett proves in \cite{paper} and \cite{book}, (p. 152 Theorem 3.13)  the following theorem:

\begin{thm}
All attractive translation invariant nearest neighbor spin systems on the integers have only the minimal and maximal invariant measures (ordered stochastically) as the extremal invariant measures. 
\end{thm}

We discuss the estimate in \cite{book} of lemma 3.7 part (e).  It is wrong, but the similar estimate $\epsilon \int g_{m,n}^{l+1}d\nu \le (4Kl+2\epsilon)\int g_{m,n}^l d\nu$ is valid.  This change does not affect the proof moving forward, because the only time the estimate is used is in the proof of lemma 3.10 in \cite{book}, where it is used to derive that $\sup_{m\le n} \int{g_{m,n}^l d\nu}<\infty$.  This fact still holds with the new estimate.  The same line of reasoning appears in \cite{paper} within the proof of lemma 2.2 there.


We will refer to \cite{book} from now on.

\section{A Correction}
The problem with the estimate as written is at the top of p.150 of \cite{book} in the discussion of bounding below the positive contribution to $\tilde{\Omega} g_{m, n}^l$.  The argument there fails to consider the possibility that change in the $\gamma$ coordinate at any of the $x_i$ among the left and right endpoints in the $l+1$ length intervals may not only create a length $l$ interval, but also destroy an adjacent length $l$ interval.  But this can only happen at at most $2g_{m, n}^l$ such sites, i.e. the left and right neighbors of intervals of length $l$. The bound below on the rate of the type of flip in question is still correctly stated as $\epsilon$ so as long as we replace $\epsilon(g_{m,n}^{l+1})$ by $\epsilon(g_{m,n}^{l+1}-2g_{m,n}^l)$ the estimate is correct.  This results in the display reading

\[\tilde{\Omega}g_{m,n}^l\ge\epsilon(g_{m,n}^{l+1}-2g_{m,n}^l)-4Klg_{m,n}^l.\]

The author would like to thank T. Liggett for his input on these matters and NSF, for part of the funding of this project.

\end{document}